\documentclass[12pt]{article}

\usepackage{amssymb,amsmath,amsthm}
\usepackage{graphicx}
\usepackage{xypic}

\theoremstyle{plain}
\newtheorem{theorem}{Theorem}
\newtheorem{prop}[theorem]{Proposition}
\newtheorem{lem}[theorem]{Lemma}
\newtheorem{cor}[theorem]{Corollary}

\theoremstyle{definition}
\newtheorem{definition}[theorem]{Definition}

\begin{document}

%
%
%
%
%
%

\newcommand{\mybox}[1]{%
	\hspace*{4pt}
  		\raisebox{-32pt}{\rule{0pt}{68pt}}
	\fbox{\makebox[28ex]{%
  		\raisebox{-24pt}{\rule{0pt}{52pt}}
		#1}}\hspace*{4pt}}

\newcommand{\eps}{\varepsilon}
\newcommand{\tensor}{\otimes}

\newcommand{\astred}{\mathop{\overline{\ast}}}
\newcommand{\astsym}{\mathop{\hat{\ast}}}

\newcommand{\Z}{\mathbb{Z}}
\newcommand{\N}{\mathbb{N}}
\newcommand{\Q}{\mathbb{Q}}

\newcommand{\WQsym}{\mathsf{WQsym}}
\newcommand{\NCQsym}{\mathsf{NCQsym}}
\newcommand{\is}{\mathsf{is}}
\newcommand{\Fgt}{\mathsf{Fgt}}
\newcommand{\Lev}{\mathsf{Inc}}
\newcommand{\st}{\mathsf{st}}
\newcommand{\Sp}{\mathsf{Sp}}
\newcommand{\Aut}{\mathsf{Aut}}
\newcommand{\Fin}{\mathsf{Fin}}
\newcommand{\Mod}{\mathsf{Mod}}
\newcommand{\spann}{\mathsf{span}}
\newcommand{\Comp}{\mathsf{Comp}}
\newcommand{\Red}{\mathsf{Red}}
\newcommand{\inv}{\mathsf{inv}}

\renewcommand{\SS}{\mathcal{S}}
\newcommand{\TT}{\mathcal{T}}


%
%
%
%
%
%

%
%
%
%
%
%
%
%
%
%
%
%
%
%
%
%
%
%
%
%

\title{\textbf{%
Trees, set compositions,\\ 
and the twisted descent algebra}}

\author{%
  \begin{minipage}{130mm}
	  \textbf{Fr\'ed\'eric Patras}\\
	  \begin{footnotesize}
		  CNRS, UMR 6621,
		  Parc Valrose,
		  06108 Nice cedex 2,
		  France\\
		  patras@math.unice.fr\\
	  \end{footnotesize}
  \end{minipage}
  \and 
  \begin{minipage}{130mm}
  \textbf{Manfred Schocker}\\
	  \begin{footnotesize}
		  Department of Mathematics,
		  University of Wales Swansea,
		  Swansea SA2 8PP,
		  UK\\
		  m.schocker@swansea.ac.uk
	  \end{footnotesize}
  \end{minipage}
  }

\date{}
\maketitle

%
%
%
%
%
%
%
%
%
%
%
%
%
%
%
%
%
%
%
%

\begin{abstract}
	\noindent
	We first show that increasing trees are in bijection with 
	set compositions, extending simultaneously a recent result 
	on trees due to Tonks and a classical result on increasing
	binary trees. We then consider algebraic structures on the 
	linear span of set compositions (the twisted descent algebra).
	Among others, a number of enveloping algebra structures are
	introduced and studied in detail. For example, it 
	is shown that the linear span of trees carries an 
	enveloping algebra structure and embeds as such in an
	enveloping algebra of increasing trees. All our constructions
	arise naturally from the general theory of twisted Hopf
	algebras. 
\end{abstract}

%
%
%
%
%
%
%
%
%
%
%
%
%
%
%
%
%
%
%
%

\section*{Introduction}
The direct sum of the Solomon-Tits algebras of type $A_n$, 
or twisted descent algebra, has been shown in~\cite{patsch06} 
to carry a rich algebraic structure which extends and generalizes
the structures on the classical descent algebra (the direct sum 
of the Solomon algebras of type $A_n$). From a combinatorial 
point of view, moving from the classical to the twisted descent
algebra means moving from the combinatorics of \emph{compositions}
(sequences of integers) to the combinatorics of 
\emph{set compositions} (sequences of mutually disjoint sets).

The purpose of the present article is to pursue further the 
study of the algebraic structures associated to set compositions, 
or twisted descents, and related objects.
Fields of application of the theory include, among others, 
the geometry of Coxeter complexes of type $A_n$, the internal
structure of twisted Hopf algebras, Markov chains associated to
hyperplane arrangements, and Barratt's twisted Lie algebra 
structures in homotopy theory. We refer to~\cite{patsch06} for 
a survey of the history of the subject, and further details on 
its various fields of application.

We first show, in Section~\ref{S1}, that the natural basis of 
the Solomon-Tits algebra of type $A_n$ (or, equivalently, 
the set of faces of the hyperplane arrangement of type $A_n$, 
the set of set compositions of $\{1,\ldots,n\}$, or the set of 
cosets of standard parabolic subgroups in the symmetric 
group $S_n$) is in bijection with the set of increasing planar
rooted trees with $n$ branchings. The result is new, to our 
best knowledge, although it appears to be a very natural 
extension of the classical bijection between permutations 
and increasing planar binary trees (see, for
example,~\cite{raw86,stan86} and the appendix of~\cite{loday01}),
and of the mapping from 
set compositions to planar rooted trees introduced 
in~\cite{tonks97} and further studied in~\cite{chapoton00}. 
This provides another link between
combinatorial structures, Hopf algebras and trees. This domain
has received a considerable attention recently, due
in particular to the discovery of 
its role in the understanding of high energy physics through 
the seminal work of Connes and Kreimer on Feynman graphs 
and Zimmermann's renormalization formula~\cite{conkre98,grafig04}.
Other approaches and other problems have also enlightened the
explanatory power of this link. These influential
contributions include Chapoton's work on Hopf algebras,
trees and the geometry of Coxeter complexes~\cite{chapoton00},
Loday's and Ronco's work on planar binary trees and
operads~\cite{lodron98}, and the work of Brouder-Frabetti
on planar binary trees and QED~\cite{broufra03}.

Set compositions of finite sets of positive integers are also
in $1$-$1$ correspondence with monomials in non-commuting
variables $x_1,x_2,\ldots$. This connects our work with the 
theory of word quasi-symmetric functions, $\WQsym$,
and the theory of quasi-symmetric functions in non-commuting
variables, $\NCQsym$, recently introduced
by Hivert et al. (see, for example,~\cite{novthibon06})
and Bergeron et al. (see, for example,~\cite{bergeron1,bergeron2}), respectively.

In Sections~\ref{S2} and~\ref{S3}, we revisit the twisted Hopf
algebra structure on the twisted descent algebra introduced in~\cite{patsch06} and start to analyze the algebraic 
implications of our combinatorial result. We study in detail 
the two Hopf algebra structures on increasing planar rooted trees
which are induced by this twisted Hopf algebra, by means of the symmetrisation and cosymmetrisation processes
of~\cite{patreu04,stover93}.
 
One of these Hopf algebras is neither commutative nor 
cocommutative. This Hopf algebra turns out to be dual to the 
algebra $\NCQsym$, thereby revealing a different approach to 
the algebra of quasi-symmetric functions in non-commuting 
variables. This observation has a striking consequence.
Namely, the classical triple at the heart of Lie theory (Solomon
algebra of type $A_n$; descent algebra (the direct sum of Solomon
algebras); quasi-symmetric functions) lifts to the world of
twisted objects as a triple (Solomon-Tits algebra of type $A_n$;
twisted descent algebra; quasi-symmetric functions in 
non-commuting variables). The link between the first two objects
was investigated in detail in~\cite{patsch06}, whereas the
duality properties between the twisted descent algebra and the
algebra of quasi-symmetric functions in non-commuting variables 
are a by-product of our considerations in Section~\ref{S2} 
of the present article. This emphasizes, once again, the 
conclusion drawn from~\cite{patreu04,patsch06,sch05}
that the twisted descent algebra is the ``natural framework''
to lift classical algebraic and combinatorial structures
(compositions, descents, shuffles, free Lie algebras, and so on) 
to the enriched setting of set compositions, tensor species,
Barratt's free twisted Lie algebras, and so on.

%


The other Hopf algebra structure provides the twisted descent
algebra (or, equivalently, the linear span of increasing planar
rooted trees), with the structure of an enveloping algebra.
This is the Hopf algebra we will be mainly interested in, 
in view of its rich algebraic and combinatorial structure. We 
show, for example,
that the twisted descent algebra is a free associative algebra 
with generators in bijection with so-called balanced increasing
rooted trees, and furthermore, naturally, the enveloping algebra 
of a free Lie algebra that we describe explicitly.


In final Section~\ref{S4}, we study Hopf subalgebras of the 
twisted descent algebra. The relationships of the noncommutative
noncocommutative Hopf algebra structure on set compositions to 
its most remarkable Hopf subalgebras appear quite simple and
natural. For example, whereas Chapoton's Hopf 
algebra~\cite{chapoton00} is related to the Malvenuto--Reutenauer 
Hopf algebra by means of a sub-quotient construction, we can 
show that the natural embedding of the free twisted associative
algebra on one generator into the twisted descent algebra
implies that the (noncommutative noncocommutative) Hopf 
algebra we consider contains Malvenuto--Reutenauer's (and 
therefore all its Hopf subalgebras) as proper Hopf subalgebras. 

However, our main concern is, once again, the other structure,
that is, the cocommutative case. We show that the direct sum of 
the symmetric
group algebras, when provided with the enveloping algebra 
structure introduced in~\cite[Section~6]{patreu04}, embeds as 
an enveloping algebra in the twisted descent algebra. We recover 
in particular, as a corollary of our results on the twisted 
descent algebra, Theorem~21 in~\cite{patreu04}, stating that the 
Lie algebra of primitive elements in this direct sum is a free 
Lie algebra. 

Finally, we also introduce enveloping algebra structures on 
planar trees and planar binary trees 
which seem to be new. We show that these 
enveloping algebras embed into the enveloping algebra of
increasing planar rooted trees as well.

%
%
%
%
%
%
%
%
%
%
%
%
%
%
%
%
%
%
%
%

\section{Set compositions and planar rooted trees}
\label{S1}

We give here a $1$-$1$ correspondence between set compositions 
and increasing planar rooted trees which extends a construction
of~\cite{tonks97} (see also~\cite{chapoton00}).
Our bijection also extends the classical correspondence between increasing planar binary rooted trees and permutations.

Let $n$ be a non-negative integer and set $[n]:=\{1,2,\ldots,n\}$.  A \emph{set composition} of $[n]$ of length $k$ is a $k$-tuple
$P=(P_1,\ldots,P_k)$ of mutually disjoint non-empty subsets
$P_1,\ldots,P_k$ of $[n]$ such that $P_1\cup\ldots\cup P_k=[n]$. 
There is an obvious $1$-$1$ correspondence between surjective maps
$\varphi:[n]\to[k]$ and set compositions $P$ of length $k$ of $[n]$
which assigns to any such $\varphi$ the $k$-tuple
$P=(\varphi^{-1}(1),\varphi^{-1}(2),\ldots,\varphi^{-1}(k))$.  For
example, if $n=3$, $k=2$ and $\varphi(1)=1=\varphi(3)$, 
$\varphi(2)=2$, then $P=(13,2)$. Note that we dropped several commas 
and curly brackets in $P$.

\begin{figure}[h]
	\hspace*{\fill}
	 $T_0=\,$\raisebox{%
	 	-4.5mm}{\includegraphics[width=12mm]{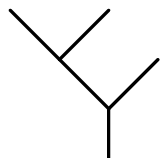}}
	\hspace*{\fill}
	 $T_1=\,$\raisebox{%
	 	-4.5mm}{\includegraphics[width=12mm]{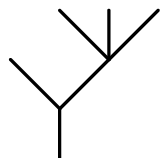}}
	\hspace*{\fill}
	 $T_2=\,$\raisebox{%
	 	-4.5mm}{\includegraphics[width=12mm]{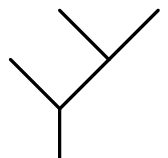}}
	\hspace*{\fill}
	\caption{\it Three trees.}
	\label{fig01}
\end{figure}

In what follows, a (planar, rooted) \emph{tree} is a finite planar
non-empty oriented connected graph $T$ without loops such that any
vertex of $T$ has at least two incoming edges and exactly one outgoing edge. In illustrations, the
root appears at the bottom, the leaves appear at the top, and the
orientation is dropped with the understanding that all
edges are oriented from top to bottom. Three trees
$T_0$, $T_1$, $T_2$ with two vertices and, respectively,
three, four and three leaves are displayed in
Figure~\ref{fig01}.
The \emph{trivial tree} 
\,\raisebox{-1pt}{%
	\includegraphics[width=1.3pt]{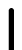}}\,
without vertex will be denoted by $\eps$.

\begin{figure}[h]
	\hspace*{\fill}
	$\bigvee(T_0,T_1,T_2)=$
	\raisebox{-10mm}{%
		\includegraphics[width=40mm]{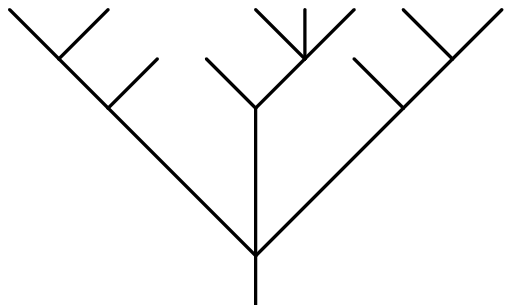}}
	\hspace*{\fill}
	\caption{\it The wedge of three trees.}
  \label{fig02}
\end{figure}

If $m\ge 1$, then the \emph{wedge} of $m+1$ trees $T_0,T_1,\ldots,T_m$ is obtained by
joining the roots of $T_0,\ldots, T_m$ to a new vertex and 
creating a new root.  It is denoted by
$\bigvee(T_0,T_1,\ldots,T_m)$. For example, the wedge of the trees 
$T_0$, $T_1$, $T_2$ in Figure~\ref{fig01}
is given in Figure~\ref{fig02}.
Any tree $T\neq\eps$ can be written uniquely as the
wedge $\bigvee(T_0,T_1,\ldots,T_m)$ of certain subtrees of~$T$. 
We define the \emph{$m$-corolla} $C_m$
by $C_m=\bigvee(\underbrace{\eps,\eps,\ldots,\eps}_{m})$. Hence $C_m$ is
the unique tree with one vertex and $m$ leaves.  For example,
$C_4=\raisebox{-3mm}{\includegraphics[width=12mm]{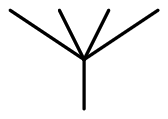}}$.

Tonks~\cite{tonks97} and Chapoton~\cite{chapoton00} studied
a surjective map which assigns to each set composition $P$ a 
tree $T$.
We will consider additional structure on trees which will allow us to turn
this surjection into a bijection.

\begin{figure}[h]
	\hspace*{\fill}
	\includegraphics[width=40mm]{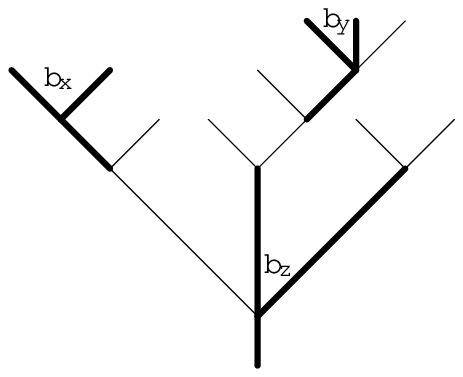}
	\hspace*{\fill}
	\caption{\it Three branchings.}
  \label{fig03}
\end{figure}

A \emph{branching} $b$ of a tree $T$ is a subgraph of $T$ 
isomorphic to 
$C_2=\raisebox{-2mm}{\includegraphics{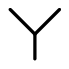}}$.
Three branchings $b_x,b_y,b_z$ are indicated in
Figure~\ref{fig03}.

\begin{figure}[h]
	\hspace*{\fill}
	$C_{m+1}=$
	\raisebox{-7mm}{%
		\includegraphics[width=50mm]{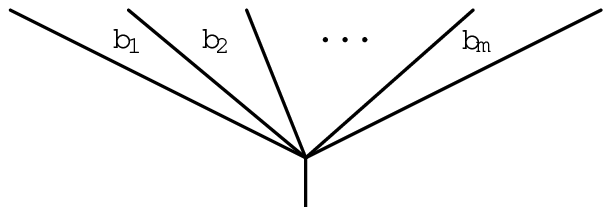}}
	\hspace*{\fill}
	\caption{\it The branchings of the $m+1$-corolla.}
  \label{fig04}
\end{figure}

We denote the set of all branchings of $T$ by $B(T)$, and its
cardinality by $b(T)$.  The set $B(T)$ admits a ``left-to-right''
order $\preceq$ in a natural way: if $T$ is the $m+1$-corolla $C_{m+1}$ for some $m\ge 0$, then $T$ has $m$ branchings $b_1,b_2,\ldots,b_m$ (labelled from left to right)
and we define $b_1\preceq b_2\preceq\cdots\preceq b_m$; 
see Figure~\ref{fig04}.
If, more generally, $T$ is the wedge of (possibly non-empty) trees 
$T_0,\ldots,T_m$,
we proceed by induction and extend the orders on
$B(T_0),\ldots,B(T_m)$ to an order on $B(T)$ by setting
$$
b_i\prec b\prec b_{i+1}
$$
for all $i\in\{0,1,\ldots,m\}$ and $b\in B(T_i)$, where $T_i$ is
embedded in $T$.  (The term $b_i$ on the left, respectively, $b_{i+1}$
on the right, does not appear when $i=0$, respectively, when $i=m$.)
For example, for the branchings indicated in Figure~\ref{fig03}, we
obtain $b_x\preceq b_y\preceq b_z$.

\begin{figure}[h]
	\hspace*{\fill}
	\includegraphics[width=40mm]{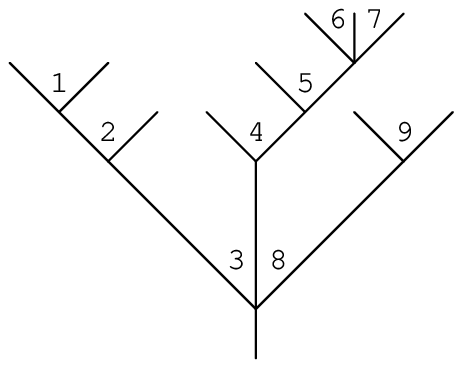}
	\hspace*{\fill}
	\caption{\it Natural labelling of the branchings.}
  \label{fig05}
\end{figure}

The order $\preceq$ induces a \emph{natural labelling} of the
branchings of $T$, namely the order-preserving map $([n],\le)\to
(B(T),\preceq)$, where $n=b(T)$ and $\le$ denotes the usual order
on $[n]$. The natural labelling of the tree in
Figure~\ref{fig03} is illustrated in Figure~\ref{fig05}.

A \emph{root branching} of $T$ is a branching $b\in B(T)$ which
contains the root. The tree in Figure~\ref{fig05}, for example, has
two root branchings, with natural labels $3$ and $8$.
With $T=\bigvee(T_0,\ldots,T_m)$ as above,
the natural label of a root branching $b_i$ is
  $$
  x_i=i+\sum_{j=0}^{i-1} b(T_j)
  $$
  for all $i\in[m]$. Furthermore, if $i\in[m]\cup\{0\}$, 
  $\tilde{b}\in B(T_i)$ has natural label $j$ in $T_i$ and
  $b$ is the branching of $T$ corresponding to $\tilde{b}$ 
  (via the embedding of $T_i$ in $T$), then $b$ has natural 
  label $x_i+j$ in $T$, where $x_0:=0$.

A \emph{level function} on a tree $T$ is a surjective map
$\lambda:V\to X$ such 
that $\lambda$ is strictly increasing along each path connecting 
a leaf of $T$ with the root of $T$. Here $V$ is the set of 
vertices of $T$ and $X$ is a totally ordered set. 
Such a level function is said to be \emph{standard}
if $X=[k]$ (with the usual order), for some non-negative 
integer $k$. A tree provided with a standard level function is 
an \emph{increasing tree}. 

\begin{figure}[h]
	\hspace*{\fill}
		\includegraphics[width=22.5mm]{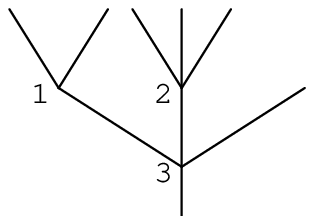}
	\hspace*{\fill}
	\hspace*{1ex}
		\includegraphics[width=22.5mm]{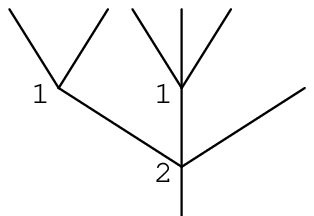}
	\hspace*{\fill}
	\hspace*{1ex}
		\includegraphics[width=22.5mm]{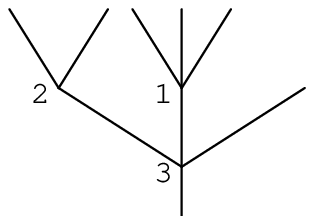}
	\hspace*{\fill}
	\caption{\it Standard level functions.}
  \label{fig06}
\end{figure}

Three standard level functions are illustrated in
Figure~\ref{fig06}.
It is easy to see that there are no other standard level 
functions on this tree.

Suppose $X$ has order $k$ and $\iota:X\to[k]$ is the
order-preserving bijection. It is clear that 
any level function $\lambda:V\to X$ on $T$
yields the standard level function $\iota\circ\lambda$
on $T$. We refer to this standard level function as 
the \emph{standardization} of $\lambda$. It will be 
advantageous at a later stage to consider arbitrary ordered 
sets $X$ in the definition of a level function.
 
\begin{figure}[h]
	\hspace*{\fill}
		\includegraphics[width=100mm]{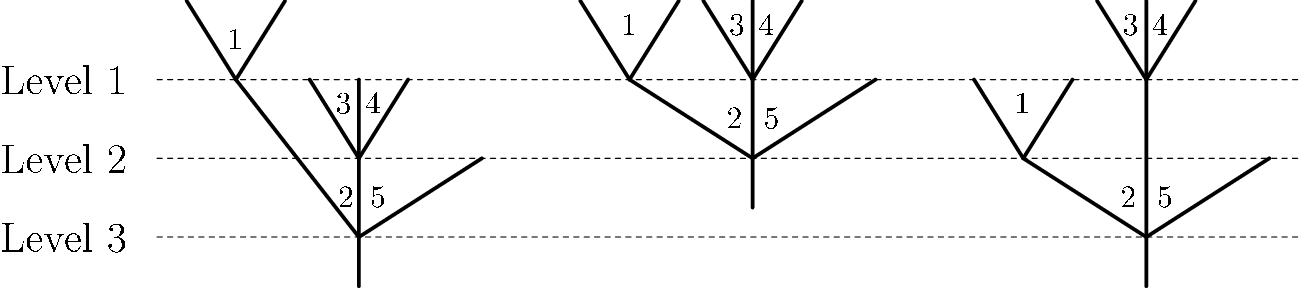}
	\hspace*{\fill}
	\caption{\it Increasing trees.}
  \label{fig07}
\end{figure}

When illustrating an increasing tree, it is more convenient to draw each vertex $v$
at the level $\lambda(v)$ (where levels increase from top to bottom)
rather than to label $v$ by $\lambda(v)$; see Figure~\ref{fig07}.
(We have added here the natural labels of the branchings in 
each case for purposes which will be clear later.)

Some geometrical properties of increasing trees can be described
in terms of natural labels and levels. For example, the vertices
of two branchings $b$ and $b'$ of an increasing tree $T$ with
natural labels $i<i'$ at levels $l<l'$ (or $l>l'$, respectively)
belong to a common path joining a leaf of $T$ to the root of $T$ 
if and only if the vertices of all the branchings $d$ with labels
between $i$ and $i'$ have levels (strictly) greater than $l$ 
(or $l'$, respectively).
 
If $\lambda:V\to X$ is a level function on $T$ and $b$ is 
a branching of $T$ with vertex $v$, we assign to $b$ the 
level $\lambda(v)$. Let
$n=b(T)$, and suppose that $X=[k]$ for some non-negative 
integer~$k$. Then the composite of the natural labelling of $B(T)$ 
with the standard level function $\lambda$ (extended to $B(T)$ in the way described) yields a
surjective map $[n]\to[k]$ or, equivalently, a set composition
$$
\sigma(T,\lambda)=P=(P_1,\ldots,P_k).
$$
More explicitly, $P_i$
consists of the natural labels of all branchings of $T$ at 
level~$i$, for all $i\in[k]$. 
The set compositions arising in this way from the increasing trees given in Figure~\ref{fig07} are, respectively,  
$(1,34,25)$, $(134,25)$ and $(34,1,25)$.

\begin{theorem} \label{bijection}
  The map $(T,\lambda)\mapsto \sigma(T,\lambda)$ is a bijection from
  the set of all increasing trees with $n$ branchings onto the 
  set of all set compositions of $[n]$.
\end{theorem}

By restriction, we obtain a correspondence between increasing
binary trees and permutations, identified with set compositions
of $[n]$ of length $n$.
Note that, in this correspondence, a permutation $\pi$
is associated to the binary tree which is classically 
associated to $\pi^{-1}$; see~\cite{raw86,stan86} 
and the appendix of~\cite{loday01}.

In order to prove the theorem, we give an inductive description 
of the inverse map $\tau$ of $\sigma$. For this purpose, it is
convenient to identify a standard level function 
$\lambda:B(T)\to[k]$
on a tree~$T$ with $n$ branchings with the corresponding map
$\varphi:[n]\to[k]$,
by means of the natural labelling of the branchings.

Let $P=(P_1,\ldots,P_k)$ be any set composition of $[n]$,
and denote the corresponding surjective map from $[n]$ to $[k]$
by $\varphi$. 
Define $T_1$ to be the 2-corolla with level function 
$\lambda_1=\varphi|_{[1]}$, the restriction of $\varphi$ 
to $[1]$. 
Let $i\in[n-1]$, and assume inductively that a tree $T_i$ has 
been constructed with level function $\lambda_i=\varphi|_{[i]}$.
Then define $(T_{i+1},\lambda_{i+1})$
to be the (unique) increasing tree 
with $i+1$ branchings obtained by adding to $(T_i,\lambda_i)$ 
a branching on the right at level $\varphi(i+1)$. 
After $n$ steps, we arrive at a tree $T_n$ with $n$ branchings 
and standard level function $\lambda_n=\varphi$. We set
$$
\tau(P)
:=
(T_n,\lambda_n).
$$
It is immediate from the definitions that $\tau$
is a left and a right inverse of $\sigma$. This proves the
theorem.

\begin{figure}[h]
	\hspace*{\fill}
	\includegraphics[width=108mm]{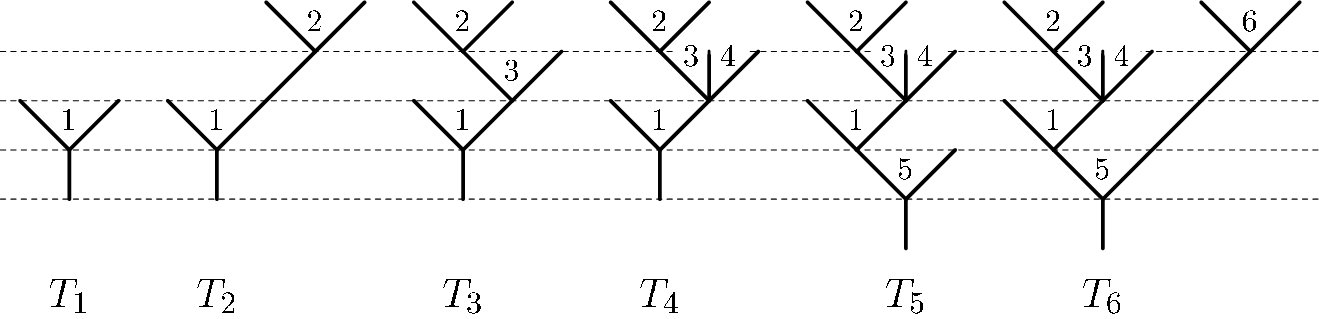}
	\hspace*{\fill}
	\caption{\it Construction of the increasing tree corresponding 
	to $(26,34,1,5)$.}
  \label{fig08}
\end{figure}

The construction is best understood through an example. 
The case where $P=(26,34,1,5)$ is illustrated in 
Figure~\ref{fig08}.

The composite of $\tau$ with the forgetful map from increasing
trees to trees can be shown to agree with Chapoton's map $\Psi$
(see~\cite[pp.267]{chapoton00}), although their definitions might look different at first sight. The proof is left to the reader.

To conclude this section, we observe that the construction 
of $\tau(P)$ applies, more generally, to an arbitrary surjective 
map $\varphi$ from a finite ordered set $A$ onto an ordered 
set $X$ (instead of the surjection $\varphi:[n]\to[k]$ 
corresponding to $P$). As a result, we obtain a tree with 
branchings labeled by $A$ and levels drawn from $X$.  
In particular, we can assign an increasing tree $\tau(Q)$ 
with labels in an ordered set $A$ to a set composition $Q$ of 
length $k$ of an arbitrary finite ordered set $A$ (since
such a set composition corresponds to a surjective 
map $\varphi:A\to[k]$).

\begin{figure}[h]
	\hspace*{\fill}
	\raisebox{-17.5mm}{%
		\includegraphics[width=35mm]{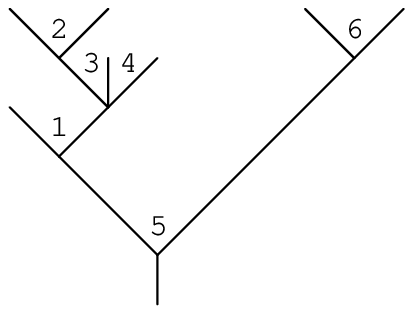}}
	\qquad
	$\longrightarrow$
	\qquad
	\raisebox{-8.75mm}{%
		\includegraphics[width=22.75mm]{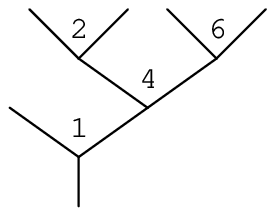}}
	\hspace*{\fill}
	\caption{\it The contraction of $(26,34,1,5)$ relative to $A=\{1,2,4,6\}$.}
  \label{fig09}
\end{figure}

If $P=(P_1,\ldots,P_k)$ is a set composition of $[n]$ and $A$ is 
a subset of $[n]$, then $Q:=(P_1\cap A,\ldots,P_k\cap A)^\#$
is a set composition of $A$, where the upper index $\#$ indicates that empty sets are deleted. In geometric terms, the increasing 
tree $\tau(Q)$ is then obtained from the increasing tree
$(T,\lambda):=\tau(P)$ by ``contracting'' in a certain way all
branchings of $T$ with natural labels not contained in~$A$ 
(and keeping the levels of branchings with labels in $A$).
Accordingly, $\tau(Q)$ is called the 
\emph{contraction of $(T,\lambda)$ relative to $A$},
or \emph{$A$-contraction of $(T,\lambda)$}. 
An example is displayed in Figure~\ref{fig09}, where
$(T,\lambda)$ is the increasing tree corresponding to
the set composition $P=(26,34,1,5)$ and $A=\{1,2,4,6\}$.
The notion of contraction will play a vital role in our
constructions of coproducts on the linear span of increasing trees
in the sections that follow.

%
%
%
%
%
%
%
%
%
%
%
%
%
%
%
%
%
%
%
%

\section{The twisted descent algebra as a Hopf algebra}
\label{S2}

Let $\Z$ and $\N$ denote the sets of all integers and of all
positive integers, respectively, and set $\N_0:=\N\cup\{0\}$.
The free twisted descent algebra $\TT$ has $\Z$-linear basis 
the set of all set compositions of finite subsets of $\N$.
Equivalently, by Theorem~\ref{bijection}, we can take the
set of all increasing trees (with natural labels drawn from $\N$)
as a basis of $\TT$.

The algebra $\TT$ is a \emph{twisted Hopf algebra}.
It was shown in~\cite{patreu04}, for arbitrary twisted
Hopf algebras $(H,\ast,\delta)$, that the twisted product
$\ast$ induces two (ordinary) products on $H$
and that the twisted coproduct $\delta$ induces two 
(ordinary) coproducts on $H$. These turn $H$ into an
ordinary Hopf algebra in two different ways.

In this section, we will make explicit these Hopf algebra
structures on $\TT$ and thereby explore the algebraic 
implications of the definitions given in~\cite{patsch06}
for the study of set compositions and increasing trees. 
In particular, we will show that one of our constructions
recovers the graded dual of the algebra of quasi-symmetric
functions in non-commuting variables considered 
in~\cite{bergeron2} (see Lemmas~\ref{3.4} and~\ref{3.7}).

For notational brevity, we refer to the free twisted descent
algebra of~\cite{patsch06} simply as ``the twisted descent 
algebra''.
All graded vector spaces considered here are connected: 
their degree~$0$ component is naturally isomorphic to the 
ground ring. Hence the two notions of bialgebra and of Hopf 
algebra coincide on these spaces, and there is no need to
specify the antipode.

Let us recall some definitions. More details and references 
on twisted algebraic structures can be found 
in~\cite{patreu04,patsch06,sch05}. 
A \emph{tensor species} is a functor from the category of finite sets and set isomorphisms $\Fin$ to the category $\Mod$ of vector spaces over a field or modules over a commutative ring. 
Unless otherwise specified, we will work over $\Z$, so that
$\Mod$ is the category of abelian groups. For convenience,
we will also assume that the finite sets we consider (and 
therefore the objects in $\Fin$) are subsets of $\N$.

The category $\Sp$ of tensor species is a linear symmetric 
monoidal category for the tensor product defined by:
$$
(F\tensor G)(S)
:=
\bigoplus\limits_{T\coprod U=S}F(T)\tensor G(U)
$$
for all $F,G\in \Sp$, $S\in \Fin$.
Here $T\coprod U=S$ means that $S$ is the \emph{disjoint} union
of $T$ and $U$.

Let $\Comp_S$ denote the set of set compositions of $S$, 
for all $S\in\Fin$. As a tensor species, the twisted descent 
algebra $\TT$ can be identified with the linearization of the 
set composition functor:
$$
\TT(S):=\Z[\Comp_S].
$$
We also write $\TT_S$ instead of $\TT(S)$.

The twisted descent algebra carries two products:
the internal or composition product $\circ$, and the
external or convolution product $\ast$. It also carries
a coproduct $\delta$. All three structures are induced by
the natural action of set compositions on twisted Hopf 
algebras~\cite{patsch06}. They are defined as follows.

\begin{definition}
	The $\Fin$-graded components $\TT_S$ of $\TT$ are associative 
	unital algebras for the \emph{composition product $\circ$},
	defined by:
	\begin{eqnarray*}
		\lefteqn{%
			(P_1,\ldots,P_k)\circ (Q_1,\ldots,Q_l)}\\[1mm]
		& := &
		(P_1\cap Q_1,\ldots,P_1\cap Q_l,
			\quad\ldots\quad,
	 	P_k\cap Q_1,\ldots,P_k\cap Q_l)^{\#}
	\end{eqnarray*}
	for all $(P_1,\ldots,P_k),(Q_1,\ldots,Q_l)\in \Comp_S$.
	As in Section~\ref{S1}, the upper index $\#$ indicates that
	empty sets are deleted. 
\end{definition}

The algebra $(\TT_{[n]},\circ)$ is now widely referred to as
the ``Solomon--Tits algebra'', a terminology introduced
in the preprint version of~\cite{patsch06}. The connections 
between Tits's seminal ideas (which ultimately led to the 
definition of the product $\circ$) and graded Hopf algebraic
structures were also first emphasized in~\cite{patsch06} and 
have since been subject to increasing interest; see, for 
example,~\cite{bergeron2} and~\cite{novthibon06}. We do not
study here the internal product and refer to~\cite{sch05}
for detailed structure results on the algebra~$(\TT_S,\circ)$.

Recall that a \emph{twisted} associative algebra is an 
associative algebra in the symmetric monoidal category of 
tensor species. Twisted versions of coassociative coalgebra,
bialgebra, and so on, are defined in the same way.

\begin{definition}
	The functor $\TT$ is provided with the structure of an 
	associative unital twisted algebra by the 
	\emph{convolution product $\ast$}, defined by:
	$$
	(P_1,\ldots,P_k)\ast (Q_1,\ldots,Q_l)
	:=
	(P_1,\ldots,P_k,Q_1,\ldots,Q_l)
	$$
	for all $S,T\in \Fin$ with $S\cap T=\emptyset$
	and all $(P_1,\ldots,P_k)\in\Comp_S$, $(Q_1,\ldots,Q_l)\in \Comp_T$.
	The identity element in $(\TT,\ast)$ is the empty 
	tuple~$\emptyset$.

	The functor $\TT$ is also provided with the structure of a 
	coassociative cocommutative counital twisted coalgebra 
	by the coproduct $\delta$, defined by:
	$$
	\delta (P_1,\ldots,P_k)
	:=
	\sum_{Q_i\coprod R_i=P_i} 
		(Q_1,\ldots,Q_k)^\# \tensor (R_1,\ldots,R_k)^\#
	$$
	for all $S\in \Fin$, $(P_1,\ldots,P_k)\in\Comp_S$.
	For example,
	\begin{eqnarray*}
		\lefteqn{%
		\delta (14,7)
		 = 
		(14,7)\tensor \emptyset
		+
		(14)\tensor (7)
		+
		(1,7)\tensor (4)}\\[1mm]
		&&
		+
		(4,7)\tensor (1)
		+
		(1)\tensor (4,7)
		+
		(4)\tensor (1,7)
		+
		(7)\tensor (14)
		+
		\emptyset\tensor (14,7).
	\end{eqnarray*}
\end{definition}

Note that
$
\delta(P)
\in 
(\TT\tensor\TT)[S]
=
\bigoplus\limits_{A\coprod B=S} \TT_A\tensor\TT_B
$
for all $P\in\Comp_S$.
We write $\delta_{A,B}$ for the component of the image 
of $\delta$ in $\TT_A\tensor\TT_B$, so that, for example, $\delta_{\{1,7\},\{4\}}(14,7)=(1,7)\tensor (4)$.

The following is an immediate consequence of~\cite{patsch06}.

\begin{prop}
	The triple $(\TT,\ast,\delta)$ is a cocommutative twisted 
	Hopf algebra. 
\end{prop}

Whenever $S,T\in\Fin$ have the same
cardinality $n$, the unique order-preserving bijection $S\to T$
induces a linear isomorphism $\is_{S,T}:\TT_S\to \TT_T$ in 
an obvious way. If, in particular, $T=[n]$, then we write 
$\TT_n:=\TT_{[n]}$ and $\is_S$ for the isomorphism from $\TT_S$ 
onto $\TT_n$.

We will now use the isomorphisms $\is_{S,T}$ to describe 
the algebra and coalgebra structures on
the twisted descent algebra which arise from the general 
constructions of~\cite{patreu04}. For this purpose, let us
introduce the graded vector space
$$
\TT_\bullet :=\bigoplus\limits_{n\in \N_0} \TT_n\,.
$$
  
\begin{lem} \label{3.4}
	The vector space $\TT_\bullet$ is a graded coassociative
	counital coalgebra with respect to the
	\emph{restricted coproduct $\overline{\delta}$},
	defined on the $n$th component by 
	$$
	\overline{\delta} 
	:=
	\bigoplus_{p+q=n}
	(\TT_p\tensor \is_{p+[q]})\circ \delta_{[p],p+[q]},
	$$
	where we write $\TT_p$ for the identity on $\TT_p$ and
	$p+[q]$ for $\{p+1,\ldots,p+q\}$. 
\end{lem}

This is a direct consequence of the definition of a 
twisted coalgebra. 

For any $A\subseteq [n]$ and any set composition
$P=(P_1,\ldots,P_k)$ of $[n]$, we write 
$$
P|_A
:=
\is_A\Big((P_1\cap A,\ldots,P_k\cap A)^{\#}\Big)
$$
so that, for example,
$
(35,62,1,47)|_{\{1,3,5,7\}}
=
\is_{\{1,3,5,7\}}((35,\emptyset ,1,7)^\# )
=
(23,1,4)
$.
Then, in particular, we have
$$
\overline{\delta}(P)
=
\sum_{p+q=n} P|_{[p]}\tensor P|_{p+[q]}
$$
for all $n\in\N_0$, $P\in\Comp_{[n]}$.  
A comparison with~\cite[Eq.~(21)]{bergeron2} now shows that
the linear map
$W:(\TT_\bullet,\overline{\delta})\to(\NCQsym^*,\Delta^*)$,
defined by $P\mapsto W_P$, is an isomorphism of coalgebras.

\begin{figure}[h]
	\begin{eqnarray*}
		\lefteqn{%
			(\TT_5\tensor \is_{5+[4]})\circ \delta_{[5],5+[4]}
			\left(
				\raisebox{	-13.3mm}{%
					\includegraphics[width=37.5mm]{figure05.eps}}
			\right)}\hspace*{38mm}\\[6mm]
		& = &
		\raisebox{-8.5mm}{%
			\includegraphics[width=33.3mm]{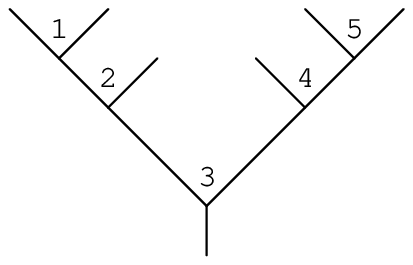}}
		\quad
		\tensor
		\quad
		\raisebox{-8.5mm}{%
			\includegraphics[width=21.6mm]{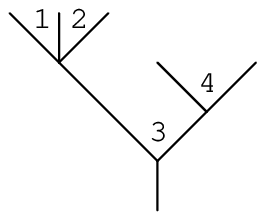}}
	\end{eqnarray*}
	\caption{\it %
		The $(5,4)$ component of the restricted coproduct
		$\overline{\delta}$.}
	\label{fig10}
\end{figure}

The restricted coproduct translates naturally into the 
language of trees, by means of Theorem~\ref{bijection}. 
If $p+q=n$, then the map 
$(\TT_p\tensor \is_{p+[q]})\circ \delta_{[p],p+[q]}$
sends an increasing tree $T$ to $T_1\tensor T_2$, where
$T_1$ is the $[p]$-contraction of $T$ and $T_2$
is the (standardization of) the $(p+[q])$-contraction of $T$,
as defined at the end of Section~\ref{S1}.
Consider, for example, the tree $T$ displayed in 
Figure~\ref{fig05}, interpreted as an increasing tree
with branchings $6$ and $7$ at level $1$, branchings $1$
and $5$ at level $2$, and so on.
The $(5,4)$ component of $\overline{\delta}(T)$ is given
in Figure~\ref{fig10}.

\begin{lem}
	The vector space $\TT_\bullet$ is a graded coassociative
	cocommutative counital coalgebra with respect to the
	\emph{cosymmetrized coproduct $\hat{\delta}$},
	defined on the $n$th component by 
	$$
	\hat{\delta} 
	:=
	\bigoplus_{A\coprod B=[n]}
	(\is_A\tensor \is_B)\circ \delta_{A,B}\,.	
	$$
\end{lem}

This follows from~\cite[pp.207]{patreu04}. Equivalently, we have
$$
\hat{\delta}(P)
=
\sum_{A\coprod B} P|_A\tensor P|_B
$$
for all $n\in\N_0$, $P\in\Comp_{[n]}$.

The cosymmetrized coproduct on increasing trees can also be
described in terms of the contraction process introduced in 
Section~\ref{S1}. If $T$ is an increasing tree with $n$ branchings
and $A\coprod B=[n]$, then the $(A,B)$ component of
$\hat{\delta}(T)$ reads $T_A\tensor T_B$, where $T_A$
and $T_B$ are the contractions of $T$ relative to $A$
and $B$, respectively.

\begin{figure}[h]
	\begin{eqnarray*}
		\lefteqn{%
			(\is_A\tensor \is_B)\circ \delta_{A,B}
			\left(
				\raisebox{-13.3mm}{%
					\includegraphics[width=37.5mm]{figure05.eps}}
			\right)}\hspace*{46mm}\\[6mm]
		& = &
		\raisebox{-11.6mm}{%
			\includegraphics[width=25mm]{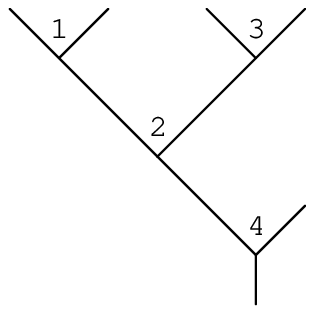}}
		\quad
		\tensor
		\quad
		\raisebox{-11.6mm}{%
			\includegraphics[width=17.5mm]{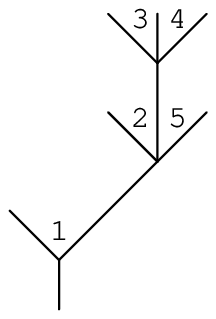}}
	\end{eqnarray*}
	\caption{\newline\it %
		The $(\{1,2,5,8\},\{3,4,6,7,9\})$ component of the 
		cosymmetrized coproduct $\hat{\delta}$.}
	\label{fig11}
\end{figure}

For example, if $T$ is the increasing tree displayed in
Figure~\ref{fig05} again and $A=\{1,2,5,8\}$, $B=\{3,4,6,7,9\}$,
then the $(A,B)$ component of $\hat{\delta}(T)$ is displayed in
Figure~\ref{fig11}.

\begin{lem} \label{lem10}
	The vector space $\TT_\bullet$ is a graded 
	associative unital algebra with respect to the 
	\emph{restricted product $\overline{\ast}$}, 
	defined on $\TT_p\tensor\TT_q$ by 
	$$
	\overline{\ast}:=\ast\circ (\TT_p\tensor \is_{[q],p+[q]})
	$$
	for all $p,q\in\N_0$.
\end{lem}

\begin{figure}[h]
	$$
	\raisebox{-6.6mm}{\includegraphics[width=13.3mm]{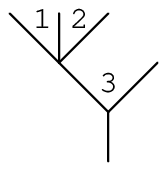}}
	\qquad
	\overline{\ast}
	\quad
	\raisebox{-6.6mm}{\includegraphics[width=29mm]{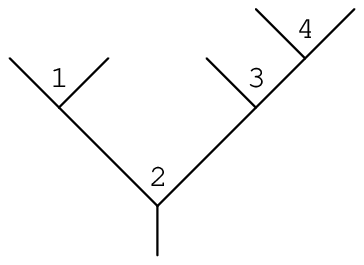}}
	\quad
	=
	\quad
	\raisebox{-6.6mm}{\includegraphics[width=41.6mm]{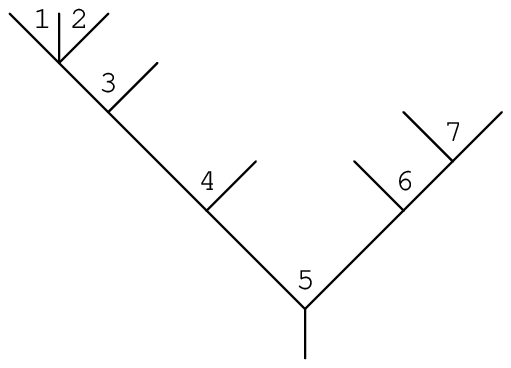}}
	$$
	\caption{\it %
		The restricted product $\overline{\ast}$.}
	\label{fig12}
\end{figure}

This is again a direct consequence of the definition of 
a twisted algebra. In terms of trees, the restricted product
is obtained by grafting an increasing tree $T_1$ with $p$ 
branchings on the left-most leaf of an increasing tree $T_2$ 
with $q$ branchings, resulting in a new tree $T$ with $p+q$
branchings. The level function on $T$ is obtained by keeping 
the levels of $T_1$ (viewed now as a subtree of $T$) and 
adding the root level $m$ of $T_1$ to all levels of $T_2$ 
(viewed also as a subtree of $T$). An example is given in 
Figure~\ref{fig12}. 

\begin{lem} \label{3.7}
	The vector space $\TT_\bullet$ is a graded 
	associative unital algebra with respect to 
	\emph{symmetrized product $\hat{\ast}$},
	defined on $\TT_p\tensor\TT_q$ by
	$$
	\hat{\ast} 
	:=
	\sum_{A\coprod B=[p+q],\,|A|=p,\,|B|=q}
		\ast\circ (\is_{[p],A}\tensor \is_{[q],B})
	$$
	for all $p,q\in\N_0$.
\end{lem}

This construction is dual to the construction of $\hat{\delta}$
(see~\cite[pp. 212]{patreu04}).
A comparison with~\cite[Eq.~(20)]{bergeron2} shows that
the map $W$ considered after Lemma~\ref{3.4} is also an 
isomorphism of algebras from $(\TT_\bullet,\hat{\ast})$
onto $\NCQsym^*$.

\begin{figure}[h]
	\begin{eqnarray*}
	\lefteqn{%
		\raisebox{-6.6mm}{%
			\includegraphics[width=13.3mm]{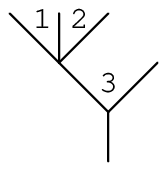}}
		\qquad
		\hat{\ast}
		\quad
		\raisebox{-6.6mm}{%
			\includegraphics[width=41mm]{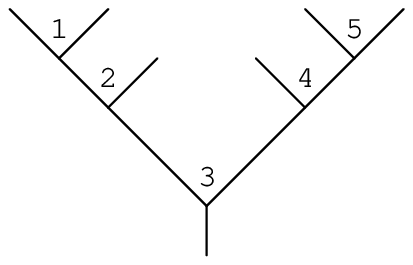}}}
		\hspace*{32mm}\\[5mm]
	& = &
	\cdots
	+
	\quad
	\raisebox{-10mm}{\includegraphics[width=48mm]{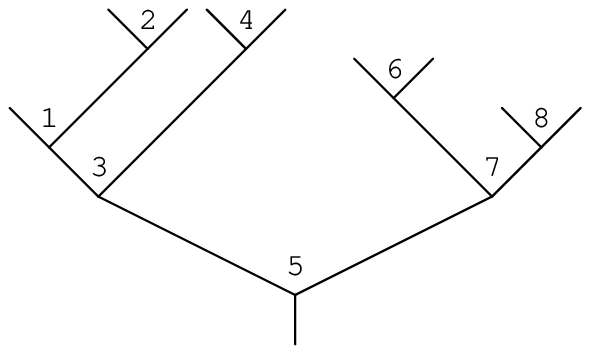}}
	\quad
	+
	\cdots
	\end{eqnarray*}
	\caption{\newline\it %
		The $(\{2,4,6\},\{1,3,5,7,8\})$ summand 
		of the symmetrized product $\hat{\ast}$.}
	\label{fig13}
\end{figure}

The symmetrized product is slightly more difficult to describe 
in terms of increasing trees and should be thought of as the 
right notion of ``shuffle product'' for increasing trees.
We give an example in Figure~\ref{fig13}.

\begin{theorem}
	The triple $(\TT_\bullet,\overline\ast ,\hat\delta )$ 
	is a graded connected cocommutative Hopf algebra.
	The triple $(\TT_\bullet,\hat\ast ,\overline\delta )$ 
	is a graded connected Hopf algebra and isomorphic to
	the graded dual of the Hopf algebra of quasi-symmetric 
	functions in non-commuting variables, $\NCQsym^*$,
	considered in~\cite{bergeron2}.
\end{theorem}

This follows from~\cite[Section~3]{patreu04} 
and the remarks after Lemmas~\ref{3.4} and~\ref{3.7}.



%
%
%
%
%
%
%
%
%
%
%
%
%
%
%
%
%
%
%
%

\section{Freeness and enveloping algebra}
\label{S3}
 
In this section, we study the algebra structures on $\TT_\bullet$
given by the reduced product $\overline{\ast}$ and the
symmetrized product $\hat{\ast}$.
Both algebras turn out to be free, and we will specify
a set of free generators.
With an eye on the Cartier-Milnor-Moore 
theorem~\cite{milmoo65} (see~\cite{pat94} for a modern 
combinatorial proof), we will then show that the cocommutative
Hopf algebra $(\TT_\bullet,\overline{\ast},\hat{\delta})$
is in fact the enveloping algebra of a free Lie algebra.

Let $n\in\N_0$ and set $\Comp_n=\Comp_{[n]}$.
We say that a set composition $(P_1,\ldots,P_k)\in\Comp_n$ 
is \emph{reduced} if there is no pair $(a,m)$ with $a<k$ 
and $m<n$ such that $\bigcup\limits_{i=1}^a P_i=[m]$. 
This corresponds to the notion of \emph{balanced tree}, 
where an increasing tree $T$ with $n$ branchings and $k$ 
levels is said to be \emph{balanced} if no pair $(a,m)$ 
with $a<k$ exists such that the branchings of $T$ with 
levels in $[a]$ are naturally labeled by the elements 
of $[m]$. 
It is clear that each set composition $P$ can be written 
uniquely as a product 
$P=P^{(1)}\overline{\ast}\cdots\overline{\ast}P^{(m)}$
of reduced set compositions 
$P^{(1)},\ldots,P^{(m)}$, where $\overline{\ast}$ is the 
restricted product defined in Lemma~\ref{lem10}.
Equivalently, each increasing tree factorizes uniquely
as a restricted product of balanced increasing trees.
This implies:

\begin{prop} \label{free-reduced}
	$(\TT_\bullet,\overline{\ast})$ is a free associative
	algebra, freely generated by the set of reduced set 
	compositions in $\Comp$.
\end{prop}

Here we write 
$
\Comp
=
\bigcup\limits_{n\in\N_0} \Comp_n
$
for the set of all set compositions of initial subsets of
$\N$.

For the next result, we work over the rational 
number field $\Q$ and consider the vector space
$\TT_\bullet^\Q=\Q\tensor_\Z \TT_\bullet$
with $\Q$-basis in bijection with $\Comp$.

\begin{cor} \label{envelope}
	The graded connected cocommutative Hopf algebra 
	$(\TT_\bullet^\Q ,\overline\ast ,\hat\delta )$ 
	is the enveloping algebra of a free Lie algebra whose 
	set of generators is naturally in bijection with the 
	set of reduced set compositions or, equivalently, of 
	balanced increasing trees. 
\end{cor}

\begin{proof}
	This follows from~\cite[Lemma~22]{patreu04} since
	all graded components of $\TT_\bullet^\Q$ have finite dimension.
\end{proof}
 
The generators of this free Lie algebra can be computed 
explicitly, using the techniques of~\cite{pat94} (see
also~\cite{pat93}). Let us write $\Red$ for the set of reduced
set compositions in $\Comp$ and $e^1$ for the logarithm of 
the identity of $\TT_\bullet$ in the convolution algebra of 
linear endomorphisms of~$\TT_\bullet$. Then the elements
$e^1(R)$, $R\in\Red$, form a set of free generators for
the primitive Lie algebra of 
$(\TT_\bullet^\Q ,\overline\ast ,\hat\delta )$. 
Further details can be found in the proof
of~\cite[Lemma~22]{patreu04}.

In concluding this section, we note that
Proposition~\ref{free-reduced} holds for the symmetrized algebra 
$(\TT_\bullet,\hat{\ast})$ as well. This was shown by 
Bergeron-Zabrocki~\cite{bergeron2} and 
Novelli-Thibon~\cite{novthibon06}, in the dual setting of 
quasi-symmetric functions in non-commuting variables.
In our approach it follows from Proposition~\ref{free-reduced}, 
by means of a standard triangularity argument: there is
a strict total ordering $\ll$ on set compositions such that, 
for all $P\in\Comp_p$, $Q\in\Comp_q$,
\begin{equation} \label{triangle}
	P\astsym Q
	\in
	P\astred Q
	+
	\spann_\Z\{\,U\in \Comp_{p+q}\,|\,P\astred Q\ll U\,\}.
\end{equation}
We give the details below, for the sake of completeness.

The order $\ll$ is defined as follows.
Consider the set $\overline{\N}$ of the positive integers 
together with the comma symbol: $\overline{\N}:=\N\cup \{,\}$. 
We extend the natural order on $\N$ to $\overline{\N}$
by putting $,<1$.
Any set composition of $[n]$ can be viewed as a word over 
the alphabet $\overline{\N}$. For example, the set composition $(145,26,3)$ translates into the word $145,26,3$. 

If $A,B\subseteq\N$, $P$ is a set composition of $A$
and $Q$ is a set composition of $B$, we set
$$
P\ll Q
$$
if $|A|<|B|$, or if $|A|=|B|$ and $P$ is smaller than $Q$ 
with respect to the lexicographic order on words over the 
alphabet $\overline{\N}$. For example, we have
$(245,169)\ll (245,178)$
because $6<7$, 
and
$(23,4569)\ll (234,1,9,5)$
because $,<4$. 

We claim that~\eqref{triangle} holds.
We start with following lemma which is clear from the definition.

\begin{lem} \label{congruence}
	Suppose $P$, $Q$, $R$, $S$ are set compositions of $A$, 
	$B$, $C$, $D$, respectively, and that 
	$A\cap  C=\emptyset=B\cap D$.
	If $|A|=|B|$, then $P\ll Q$ implies that $P\ast R\ll Q\ast S$.
\end{lem}

To see~\eqref{triangle} now, note first that
$P\astsym Q-P\astred Q$
is equal to the sum of all set compositions
$$
U=\is_{[p],A}(P)\ast\is_{[q],B}(Q),
$$
where $A$, $B$ are such that
$A\coprod B=[p+q]$, $|A|=p$ and $A\neq[p]$.
Hence, by Lemma~\ref{congruence}, it suffices to show that 
\begin{equation} \label{helper}
	P\ll \is_{[p],A}(P)
\end{equation}
for all subsets $A$ of $[p+q]$ of order $p$ with $A\neq[p]$.
Suppose $P=(P_1,\ldots,P_k)$, and let $\iota:[p]\to A$ be the 
order-preserving bijection. Then
$\is_{[p],A}(P)=(\iota(P_1),\ldots,\iota(P_k))$.
Hence, if $i\in[k]$ is minimal with $P_i\neq\iota(P_i)$
and $j\in P_i$ is minimal with $j\neq\iota(j)$,
then $j<\iota(j)$ because $\iota$ preserves the orders
on $[p]$ and $A$. This implies~\eqref{helper}, hence
also~\eqref{triangle}.


\begin{theorem}  \label{free-sym}
	$(\TT_\bullet,\hat{\ast})$ is a free associative
	algebra, freely generated by the set of reduced set 
	compositions in $\Comp$.
\end{theorem}

\begin{proof}
	We need to show that the symmetrized products
	$P^{(1)}\astsym\cdots\astsym P^{(m)}$
	of reduced set compositions $P^{(1)},\ldots,P^{(m)}$
	in $\Comp$ form a $\Z$-basis of~$\TT_\bullet$.
	From Lemma~\ref{congruence} and~\eqref{triangle}
	it follows that such a product is contained in
	$$
	P^{(1)}\astred\cdots\astred P^{(m)}
	+
	\spann_\Z\{%
		\,U\in \Comp\,|\,P^{(1)}\astred\cdots\astred P^{(m)}\ll U\,\}.
	$$
	Hence the claim follows from Proposition~\ref{free-reduced},
	because the reduced products 
	$P^{(1)}\astred\cdots\astred P^{(m)}$
	form a $\Z$-basis of $\TT_\bullet$.
\end{proof}

%
%
%
%
%
%
%
%
%
%
%
%
%
%
%
%
%
%
%
%

\section{Enveloping algebras and trees}
\label{S4}

In this section we consider a number of subspaces of 
the twisted descent algebra (up to natural identifications),
such as the direct sum 
$\SS_\bullet :=\bigoplus\limits_{n\in\N_0}\Z[S_n]$ of the symmetric group algebras or the linear span of (planar rooted)
trees. We will show how the 
Hopf algebra structures on set compositions and increasing trees 
studied in the previous sections restrict to these subspaces. 
In this way, we will obtain various algebraic structures.
Some of them are known: the Malvenuto-Reutenauer Hopf algebra of permutations, and the enveloping algebra structure on $\SS_\bullet$ introduced in~\cite{patreu04}; 
others seem to be new.
As far as we can say, there is no obvious connection between 
our Hopf algebras of trees
and
the Hopf algebra structures on trees and forests appearing in renormalization theory (see, for
example,~\cite{broufra03,conkre98,grafig04}).

First of all, we show that the enveloping algebra 
with underlying graded vector space $\SS_\bullet$,
as defined in~\cite{patreu04}, is a 
sub-enveloping algebra of the enveloping algebra of twisted 
descents 
$(\TT_\bullet,\overline{\ast},\hat{\delta})$
considered in Corollary~\ref{envelope} (up to an anti-involution).
Second, the Malvenuto-Reutenauer Hopf algebra is a 
Hopf subalgebra of the (non-cocommutative) Hopf algebra of
twisted descents $(\TT_\bullet,\hat{\ast},\overline{\delta})$.

Third, we turn to (planar rooted) trees and show that the corresponding graded vector space can be provided with an enveloping algebra structure. It embeds, as an enveloping algebra, into the enveloping algebra of increasing trees. 

Finally, we show that a stronger result holds for (planar rooted) binary trees, whose associated sub-enveloping algebra of the enveloping algebra of trees also embeds into the enveloping algebra associated to~$\SS_\bullet$.

\begin{figure}[h]
	\hspace*{-8pt}
	$
	\xymatrix@C=8ex@R=8ex{%
	  \mybox{Planar rooted trees}
		\ar@{^{(}->}[r]%
	  &
	  \mybox{%
		  \begin{minipage}{29ex}
				\begin{center}
					Increasing\\
					planar rooted trees\\ 
					=\\ 
					set compositions
				\end{center}	
			\end{minipage}}%
		\ \\ 
	  \mybox{Planar rooted binary trees}
		\ar@{^{(}->}[u]%
		\ar@{^{(}->}[r]%
	  &
	  \mybox{Permutations}
		\ar@{^{(}->}[u]%
	}
	$
	\caption{\it Diagram of embeddings.}
	\label{diagram}
\end{figure}
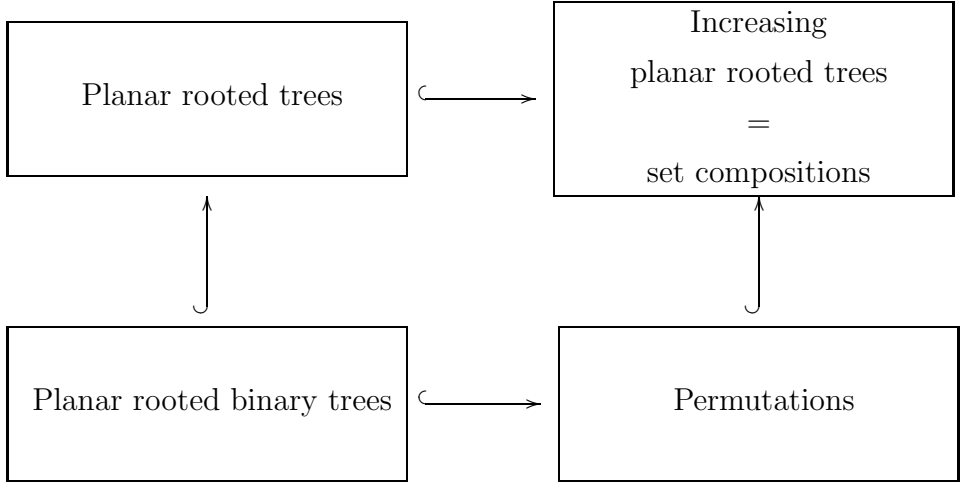

These various embeddings of enveloping algebras are 
induced by embeddings of sets of combinatorial objects (trees, 
set compositions, and so on) and illustrated in
Figure~\ref{diagram}.
Note that the embeddings are not canonical; they will be defined
below.

As far as the space $\SS_\bullet$ is concerned, all our results build on the functorial superstructure associated to the
corresponding tensor species~$\SS$. 
This tensor species is defined on objects by  
$$
\SS[T]:=\Z[\Aut_{\Fin}(T)]
$$
for all $T\in\Fin$. If $T=\{t_1,\ldots,t_k\}\subseteq \N$ 
such that $t_1<\cdots<t_k$, we can think of any
bijection $\sigma\in\Aut_\Fin(T)$ as the  
ordered sequence $(\sigma (t_1),\ldots,\sigma (t_k))$. 
The action of a bijection $\phi$ from $T$ to $S$ is then given 
by
$$
\phi (\sigma (t_1),\ldots,\sigma (t_k))
:=
(\phi\circ \sigma (t_1),\ldots,\phi\circ\sigma (t_k)).
$$
A twisted Hopf algebra structure on $\SS$ can be defined 
as follows.
If
$S=\{s_1,\ldots,s_k\},T=\{t_1,\ldots,t_l\}\subseteq\N$
such that $S\cap T=\emptyset$
and
$\alpha\in \Aut_{\Fin}(S)$, $\beta\in\Aut_{\Fin}(T)$,
then 
$$
(\alpha (s_1),\ldots,\alpha (s_k))
\times 
(\beta (t_1),\ldots,\beta (t_l))
:=
(\alpha (s_1),\ldots,\alpha (s_k),\beta (t_1),\ldots,\beta (t_l)).
$$
Furthermore, if $\sigma\in\Aut_\Fin(S\coprod T)$, then
$$
\delta_{S,T}(\sigma):=\sigma|_S\tensor\sigma|_T\,.
$$
Here we write $\sigma|_S$ for the subsequence of $\sigma$
associated to the elements of $S$. For example,
$(3,5,2,4,1)|_{\{1,3,5\}}=(3,5,1)$.

The following lemma is a direct consequence of the definitions.

\begin{lem} \label{twisted-emb}
	The canonical embeddings of the symmetric groups $S_n$ into 
	$\Comp_n$, defined by
	\begin{equation} \label{perm-emb}
		\sigma
		\longmapsto 
		(\sigma (1),\ldots,\sigma (n))
	\end{equation}
	for all $n\in\N_0$ and $\sigma\in S_n$, yield an embedding 
	of the twisted Hopf algebra $\SS$ into the twisted Hopf 
	algebra $\TT$ of set compositions.
\end{lem}

Recall that to each twisted Hopf algebra is associated a
symmetrized Hopf algebra and a cosymmetrized Hopf algebra.
The latter is an enveloping algebra (that is, a graded connected cocommutative Hopf algebra) if the twisted Hopf algebra is cocommutative.
Since
these constructions are natural in a functorial sense, an embedding of twisted Hopf algebras 
induces embeddings of the associated symmetrized and 
cosymmetrized Hopf algebras.

The cosymmetrized Hopf algebra stucture on $\SS_\bullet$ 
associated with the twisted Hopf algebra $\SS$ 
is defined as follows.

\begin{definition}
	The graded vector space $\SS_\bullet$ is an enveloping algebra with respect
	to the usual concatenation product $\times$ and
	the cosymmetrized coproduct $\hat\delta$.
	These are defined by
	$$
	(\alpha\times\beta )(i)
	:= 
	\left\{\begin{array}{ll}
		\alpha(i),     & \mbox{ if $i\le n$,}\\[1mm]
		n+\beta (i-n), & \mbox{ if  $i>n$,}
	\end{array}\right.
	$$
	for all $n,m\in\N_0$, $\alpha\in S_n$, 
	$\beta\in S_m$, $i\in[n+m]$,
	and
	$$
	\hat{\delta}(\alpha)
	:=
	\sum\limits_{S\coprod T=[n]}
	\is_S(\alpha|_S)\tensor \is_T(\alpha|_T)
	$$
	for all $n\in\N_0$, $\alpha\in S_n$.
\end{definition}  

Here $\is_S$ and $\is_T$ are the standardization maps considered
earlier so 
that, for example, $\is_{\{1,3,5\}}((3,5,2,4,1)|_{\{1,3,5\}})=(2,3,1)$.

A detailed description of the Lie algebra of primitive elements 
associated with this enveloping algebra is given 
in~\cite{patreu04}.

The symmetrized Hopf algebra stucture on $\SS_\bullet$ 
associated with the twisted Hopf algebra $\SS$ 
yields the Malvenuto-Reutenauer algebra. We recall
its definition.

\begin{definition}
	The \emph{Malvenuto-Reutenauer algebra} $\SS_\bullet$ is a 
	graded connected Hopf algebra with respect to the 
	convolution product $\ast$ and the restricted coproduct
	$\overline{\delta}$. These are defined by
	$$
	\alpha\ast\beta 
	:= 
	q_{(n,m)}\cdot (\alpha\times \beta )
	$$
	for all $n,m\in\N_0$, $\alpha\in S_n$, $\beta\in S_m$,
	and 
	$$
	\overline{\delta}(\alpha)
	:=
	\sum\limits_{i=0}^n 
	\alpha|_{[i]}
	\tensor 
	\is_{\{i+1,\ldots,n\}}(\alpha|_{\{i+1,\ldots,n\}})
	$$
	for all $n\in\N_0$, $\alpha\in S_n$.	
\end{definition}

Here we write $q_{(n,m)}$ for the sum in $\Z[S_{n+m}]$
of all permutations $\pi\in S_{n+m}$ such that
$\pi(1)<\cdots<\pi(n)$ and $\pi(n+1)<\cdots<\pi(n+m)$.

Note that these definitions agree with the structures studied
in~\cite{patreu04} only up to the involution 
$\inv:\SS_\bullet\to\SS_\bullet$ 
which maps any permutation to its inverse; 
for, in that article, twisted bialgebras were studied from 
Barratt's point of view~\cite{barratt77} (that is, by 
considering \emph{right} modules over symmetric groups or,
equivalently, by considering \emph{contravariant} functors 
from the category of finite sets and bijections), whereas 
here and in~\cite{patsch06}, twisted bialgebras have been 
studied from Joyal's point of view~\cite{joyal86} 
(that is, by considering
\emph{left} modules over symmetric groups or, equivalently, by
considering \emph{covariant} functors). As usual, one can move from one 
point of view to the other using the map $\inv$. 
Details on the two point of views and their relative behaviours 
can be found in the first section of~\cite{patreu04}.  

From Lemma~\ref{twisted-emb}, we can now deduce without further
ado:

\begin{theorem} \label{ord-emb}
	The canonical embedding of $\SS_\bullet$ into $\TT_\bullet$
	given by~\eqref{perm-emb} is an embedding of enveloping
	algebras 
	$
	(\SS_\bullet,\overline{\ast},\hat{\delta})
	\to
	(\TT_\bullet,\overline{\ast},\hat{\delta})
	$.
	It is also an embedding of the Malvenuto-Reutenauer Hopf 
	algebra into $(\TT_\bullet,\hat{\ast},\overline{\delta})$.
\end{theorem}

We now turn to the enveloping algebra structures on trees and 
binary trees and consider $\TT_\bullet$ as the linear span
of all increasing trees. Recall that the product of two increasing 
trees in the twisted algebra of increasing trees is the grafting 
of the first tree on the left most leaf of the second. If we 
assume, for simplicity, that the first tree has levels 
$1,\ldots,n$ and the second tree has non-standard levels $n+1,\ldots,n+m$, 
then the level of each branching in the product is the same
before and after the grafting.

The same operation (grafting on the left most 
leaf) defines an associative product on the linear span of 
(planar rooted) trees. The forgetful map 
$$
\Fgt:\TT_\bullet\to \overline{\TT}_\bullet
$$
from increasing trees to trees is clearly an algebra map,
where we write  
$$
\overline{\TT}_\bullet 
=
\bigoplus\limits_{n\in\N_0}\overline{\TT}_n
$$
for the graded vector space with basis the set of trees, 
graded by the number of branchings.

We claim that $\Fgt$ has a section $\Lev$ in the category of graded
associative algebras with identity.
Let $T$ be any non-empty planar rooted tree with $n$ branchings.
Then $T$ can be written uniquely as a wedge 
$T=\bigvee (T_0,\ldots,T_m)$. 
Let $b_1,\ldots,b_m$ denote the natural labels of the root
branchings of $T$. We assume inductively that an increasing tree
$(T_i,\lambda_i)=\Lev(T_i)$
has been defined for all $0\le i\le m$.
A level function $\lambda$ on $T$ can then be defined by
requiring that: 
\begin{itemize}
	\item[(i)] 
		$(T_i,\lambda_i)$ is the $S_i$-contraction of
		$(T,\lambda)$ for all $0\leq i\leq m$, 
		where $S_i$ is the set of branchings 
		of $T$ which belong to $T_i$ (embedded in $T$).
	\item[(ii)]
		The levels on $T_i$ (embedded in $T$) are strictly 
		less than the levels on $T_j$ (embedded in $T$), for 
		all $0\leq i<j\leq m$.
\end{itemize}
We set $\Lev(T):=(T,\lambda)$ and observe:

\begin{lem} \label{lev}
	The map $\Lev$ from trees to increasing trees is a section 
	of the forgetful map $\Fgt$. It defines an embedding of 
	algebras 
	$
	(\overline{\TT}_\bullet,\overline{\ast})
	\to 
	(\TT_\bullet,\overline{\ast})
	$
	where $\overline{\ast}$ is the left grafting product 
	on $\overline{\TT}_\bullet$ (by slight abuse of notation)
	and the restricted product on $\TT_\bullet$.
\end{lem}

The proof is geometrically straightforward and left to the reader.

Any increasing tree in the image of $\Lev$ is called 
\emph{left increasing}. Due to the recursive definition of $\Lev$, 
left increasing trees can be characterized as follows.

\begin{lem} \label{inc-char}
	Let $(T,\lambda)$ be an increasing tree.	Then $(T,\lambda)$ 
	is left increasing if and only if, for any branchings $b$ 
	and $b'$ of $T$, we have $\lambda(b)<\lambda(b')$
	whenever $b$ is to the left of $b'$ in $T$ and the associated
	vertices $v$, $v'$ do not lie on a common path connecting a 
	leaf of $T$ with the root of $T$.
\end{lem}

\begin{cor}
	Any contraction of a left increasing tree is left increasing.
\end{cor}

\begin{proof}
	It is enough to understand how all three notions occuring 
	in Lemma~\ref{inc-char}
	(to lie on a common path, to lie further to the left, to have a 
	smaller level) behave with respect to the contraction process.
	
	Let $(T',\lambda')$ be a contraction of a left increasing 
	tree $(T,\lambda)$.	
	Let $b$, $b'$ be branchings of $T'$, and denote the corresponding
	branchings of $T$ by $\tilde{b}$ and $\tilde{b}'$, respectively.
	
	The contraction process is an order preserving map with 
	respect to the natural labelling of branchings. This
	follows directly from its recursive left-to-right definition.
	In particular, 
	$b$ is to the left of $b'$ in $T'$ if and only if $\tilde{b}$ 
	is to the left of $\tilde{b}'$ in $T$.

	Furthermore, we have $\lambda'(b)<\lambda'(b')$ if and 
	only if $\lambda(\tilde{b})<\lambda(\tilde{b}')$ since
	$\lambda'$ is the standardization of the restriction 
	of $\lambda$ to the set of branchings of $T'$.
	 
	Let us assume that $b$ is to the left of $b'$. As mentioned
	in Section~\ref{S1}, the vertices associated to $\tilde{b}$ 
	and $\tilde{b}'$ lie on a common path from a leaf to the root 
	of $T$ if and only if either
	$\lambda(\tilde{b})<\lambda(\tilde{b}')$
	and for any branching $k$ between $\tilde{b}$ and $\tilde{b}'$ 
	in the 
	left-to-right ordering $\lambda(k)>\lambda(\tilde{b})$, or
	$\lambda(\tilde{b})>\lambda(\tilde{b}')$ and, with the same
	notation, $\lambda(k)>\lambda(\tilde{b}')$. Since the same
	characterization holds for $b$ and $b'$ in $T'$, and since levels 
	(up to standardization) and the left-to-right ordering are
	preserved by the contraction process, it follows that the
	property of lying on a common path from a leaf to the root is
	preserved by the contraction $T\mapsto T'$. 

	In particular, if $b$ and $b'$ do not lie on a common path, 
	the same property is true for $\tilde{b}$ and $\tilde{b}'$. 
	Since $T$ is left increasing, we get
	$\lambda(\tilde{b})<\lambda(\tilde{b}')$ 
	and 
	$\lambda(b)<\lambda(b')$, 
	which concludes the proof.
\end{proof}

As a consequence of the preceding result, 
the cosymmetrized coproduct $\hat{\delta}$
restricts to a coproduct on the $\Z$-linear span of left
increasing trees, $\Lev(\overline{\TT}_\bullet)$, in $\TT_\bullet$.
Combined with Lemma~\ref{lev}, this gives:

\begin{theorem} \label{inc-tree}
	The linear span of left increasing trees (or, equivalently, 
	the linear span of trees) is a Hopf subalgebra of the 
	cocommutative Hopf algebra of increasing trees. 
	In particular, this Hopf algebra of trees is an enveloping
	algebra. It is free as an algebra, and the enveloping algebra 
	of a free Lie algebra. 
\end{theorem}

To conclude, we observe that the left increasing tree 
$\Lev(T)$ corresponding to a \emph{binary} tree $T$ is 
characterized by the property that it has a single grafting 
at each level. 
This property is also preserved by the contraction process.
Furthermore, the set composition $(\sigma(1),\ldots,\sigma(n))$
corresponding to a permutation $\sigma\in S_n$ can be
characterized in the same way.
Hence we get from Theorem~\ref{inc-tree}:

\begin{cor}
	The linear span of (planar rooted) binary trees is 
	naturally embedded in $\SS_\bullet$ and $\TT_\bullet$ 
	as an enveloping algebra, and is the enveloping algebra 
	of a free Lie algebra. 
\end{cor}

As far as we can say, there is no direct connection between
our Hopf algebra of planar binary trees and 
Loday--Ronco's~\cite{lodron98}.
 
%
%
%
%
%
%
%
%
%
%
%
%
%
%
%
%
%
%

\end{document}